\documentclass{amsart}

\usepackage{hyperref, color}

\newtheorem{theorem}{Theorem}[section]

\theoremstyle{remark}
\newtheorem{remark}[theorem]{Remark}

\newcommand{\field}[1]{\mathbb{#1}}
\newcommand{\real}{\field{R}}

\newcommand{\meta}[2]{\langle #1,#2 \rangle }
\newcommand{\n} {{\bf n}}

       \newcommand{\De}{\Delta}

\begin{document}

\title[Uniqueness Theorems for fully nonlinear equations]
{Uniqueness Theorems for fully nonlinear conformal equations on subdomains of the sphere}


\author[Cavalcante]{M. P. Cavalcante}
\address{Instituto de Matem\'atica - Universidade Federal de Alagoas, Macei\'o - Brazil}
\email{marcos@pos.mat.ufal.br}
\thanks{The authors were partially supported  by CNPq-Brazil. The second author is partially supported by Spanish MEC-FEDER Grant MTM2013-43970-P.}

 \author[Espinar]{J.  M. Espinar}
 \address{Instituto Nacional de Matem\'atica Pura e Aplicada, Rio de Janeiro - Brazil}
 \email{jespinar@impa}


\date{\today}



\begin{abstract}
{In this paper we prove classification results to elliptic fully nonlinear conformal equations on certain subdomains of the sphere with prescribed constant mean curvature on its boundary. Such subdomains are the hemisphere (or a geodesic ball on $\mathbb{S}^n$) of dimension $n\geq 2$ with prescribed constant mean curvature on its boundary, and annular domains with minimal boundary.}
Our results extend the classifications of Escobar  in \cite{E0} when $n\geq 3$, and Hang-Wang in \cite{HaWa} and Jimenez in \cite{J} when $n=2$.

\end{abstract}

\maketitle

\section{Introduction}\label{intro}

Let $(M^n, g_0)$ be a compact orientable Riemannian manifold with smooth boundary and dimension $n\geq 3$. Let us denote by  $R(g_0)$  its scalar curvature and by $h(g_0)$ its boundary mean curvature with respect to the outward unit normal vector field. If $g=u^{\frac{4}{n-2}}g_0$ is a  metric  conformal to $g_0$ then its scalar curvature and
boundary mean curvature  are related by the following nonlinear elliptic partial differential equation of  critical Sobolev exponent in terms of the positive function $u$  
\begin{equation}\label{pde}
   \left\{ \begin{array}{lrc}
  \De_{g_0} u - \frac{n-2}{4(n-1)}R(g_0) u +
  \frac{n-2}{4(n-1)}R(g)u^{\frac{n+2}{n-2}}=0 &  \textrm{ in } M, \\ \\
  \frac{\partial u}{\partial \eta}+\frac{n-2}{2}h(g_0) u -
  \frac{n-2}{2}h(g)u^{\frac{n}{n-2}}=0 & \textrm{ on } \partial M,
   \end{array}  \right.
  \end{equation}
where $\Delta_{g_0}$ is the Laplace operator with respect to the metric $g_0$ and $\eta$ is the outward unit normal vector field to $\partial M$. The problem of existence of solutions to equation (\ref{pde}) when $R(g)$ and $h(g)$ are constants is referred as the \index{Yamabe problem}\emph{Yamabe problem} which was completely solved when $\partial M = \emptyset$ in a sequence of works, beginning with H. Yamabe himself ~\cite{Y}, followed by N. Trudinger
~\cite{T} and T. Aubin ~\cite{A}, and finally by R. Schoen ~\cite{S2}. 

In the case of nonempty  boundary  almost all the cases were solved by the works of J. Escobar  ~\cite{E1} and  ~\cite{E2},  Z. Han and Y. Li ~\cite{HL} and F. Marques ~\cite{M} among others. The problem of existence of solutions to 
(\ref{pde}) when $M$ is the round sphere $\mathbb S^n$ and  $R(g)$ is a prescribed  function is referred as the \emph{Nirenberg problem} and it is still open in all its generality.

The most important case to be considered in problem (\ref{pde}) is when $M$ is the unit  Euclidean ball or equivalently,  
the closed hemisphere $\mathbb S^n_+$ endowed with the standard round metric $g_0$ and $R(g)$ is a positive constant. Using the conformal invariance of the problem  we may consider only the case $R(g)=1$. Also, for simplicity we will state this problem in terms of the eigenvalues of the Schouten tensor.

Namely, given a Riemannian manifold $(M ^n , g)$ with $n\geq 3$, the Riemann curvature tensor can be decomposed as
$$ {\rm Riem}_g = W_g + {\rm Sch}_g \odot g , $$where $W_g$ is the Weyl tensor,
$\odot$ is the Kulkarni-Nomizu product, and
\begin{equation}\label{schouten}
{\rm Sch}_g := \frac{1}{n-2}\left( {\rm Ric}_g - \frac{R(g)}{2(n-1)}g\right) 
\end{equation}
is the Schouten tensor. Here ${\rm Ric}_g$ stands for the Ricci curvature of $g$. The eigenvalues of ${\rm Sch}_g$ are defined as the eigenvalues of the endomorphism $g^{-1}{\rm Sch}_g$. Note that ${\rm Trace}(g^{-1}{\rm Sch}_g) = 2 (n-1) R(g)$.

Now, take $g=e^{2\rho}g_0$  a conformal metric and denote by $\lambda (p)=(\lambda_1(p),\ldots, \lambda_n (p))$ the eigenvalues of the Schouten  tensor  of $g$ at $p \in \mathbb{S}^n$. The Yamabe problem for $R(g)=1$  and $h(g)=c$ is equivalent to find a smooth function $\rho$ on $\mathbb S^n_+$ such that

\begin{equation}\label{yamabe2}
\left\{ \begin{array}{lrc}
\frac{1}{2(n-1)} \left( \lambda_1 + \dots + \lambda_n \right) =1&\textrm{in } \mathbb S^n_+, \\ \\
h(g)=c & \textrm{on }\partial \mathbb{S}^n_+.
\end{array}\right.
\end{equation}

In \cite{E0}, Escobar showed that the solutions of problem (\ref{yamabe2}) are given by  $g = \Phi ^* g_0$ on $\overline{\mathbb{S}^n_+}$, where $\Phi$ is a conformal diffeomorphism of $(\mathbb S ^n, g_0)$.

Posed in this form, problem (\ref{yamabe2}) can be generalized for other functions of the eigenvalues of the Schouten tensor. For instance, one may consider the $\sigma_k$-Yamabe problem on $\mathbb S^n_+$  considering the 
$k$-symmetric function of the eigenvalues of the Schouten tensor. In this paper we are interested in the fully nonlinear case of this problem, in the line opened by A. Li and Y.Y. Li \cite{LL}. Namely, given $(f,\Gamma)$  an elliptic data and 
$c\in \real$, find $\rho \in C^\infty (\overline {\mathbb{S}^n_+})$
so that $g=e^{2\rho}g_0 $ is a solution of the problem

\begin{equation}\label{p1}
\left\{\begin{array}{rl}
f(\lambda(p))= 1,& \lambda(p)\in\Gamma, \, \, p\in \mathbb{S}^n_+,\\ \\

h(g) = c	, & \textrm{ on } \partial \mathbb{S}^n_+ .
\end{array} \right. 
\end{equation}

First, we will show

\begin{theorem}\label{thp1}
Let $\rho \in C^\infty (\overline{\mathbb S ^n _+})$ be a solution of \eqref{p1}. 
Then,  $g= e^{2\rho}g_0$ is given by $\Phi ^* g_0$ on 
$\overline{\mathbb{S}^n_+}$, where $\Phi$ is a conformal diffeomorphism of $(\mathbb S ^n , g_0)$.
\end{theorem}

Let us explain the meaning of a symmetric function $f$ being elliptic. Denote
$$ \Gamma _n = \{ (x_1, \ldots , x_n) \in \real ^n \, : \,\, x_i > 0, \, i=1, \ldots , n\} $$and
$$\Gamma _1 = \{ (x_1, \ldots , x_n) \in \real ^n \, : \,\, \sum _{i=1}^nx_i > 0, \, i=1, \ldots , n\} .$$
Let $\Gamma$ be an open connected component of 
$$  \{ (x_1, \ldots , x_n) \in \real ^n \, : \,\, f(x_1, \ldots , x_n) > 0\} . $$

We say that $(f,\Gamma)$ is \emph{elliptic} if they satisfy  

\begin{enumerate}

\item[(i)] $\Gamma_n\subset \Gamma\subset \Gamma_1$; \\

\item[(ii)] For all $(x_1, \ldots , x_n)\in \Gamma$ and all $(y_1, \ldots , y_n) 
\in \Gamma  \cap ((x_1,\ldots ,x_n) + \Gamma _n)$, there exists a curve 
$\gamma$ connecting $(x_1, \ldots , x_n)$ to $(y_1, \ldots , y_n)$ inside 
$\Gamma$ such that $\gamma ' \in \Gamma _n $ along $\gamma$; \\

\item[(iii)] $\dfrac{\partial f}{\partial x_i}>0$, \quad $\forall \, i=1\ldots, n;$

\item[(iv)] There exists $\lambda _0>0$ such that $f(\lambda _0, \ldots , \lambda_0 ) =1$.
\end{enumerate}

\medskip

We point out that we do not ask any convexity on $f$.  One important example that satisfies the above conditions is when we consider $f = \sigma _k$, being $\sigma _k$ the $k-$th elementary symmetric polynomial (cf. \cite{CHY} and references therein). Moreover, we require $\rho$ be $C^\infty$ to apply, avoiding technicalities, the Maximum Principle. Certainly, the condition on the regularity on $\rho$  can be relaxed, but this analytic issue is not relevant for the main ideas of the paper.  The important point is that, under these conditions, we can apply the Maximum Principle (cf. \cite{LL}).

\medskip

\begin{remark}
{\it Condition (iv) means that the Schouten tensor of a dilation of the standard metric of the sphere satisfies the equation.}
\end{remark}

We will prove Theorem \ref{thp1} using a geometric method developed by the second author, G\'alvez and Mira in 
\cite{EGM}  and the idea goes as follows.  Given a conformal metric on $\overline{\mathbb S ^n _+}$, we will construct
a horospherically concave hypersurface $\Sigma \subset \mathbb H ^{n+1}$ with compact boundary $\partial \Sigma$. By construction $\Sigma $ and $\partial \Sigma$ are topologically  $\mathbb S ^n _+$ and $\partial \mathbb S^n_+ = \mathbb S ^{n-1}$ respectively.

We next show that the boundary $\partial \Sigma$ is contained on an equidistant hypersurface to a totally geodesic hyperplane and makes a constant angle with this equidistant hypersurface along $\partial \Sigma$. Moreover, $\Sigma$ is completely contained in one of the half-spaces determined by the equidistant. This fact, jointly with a convexity argument, will prove that $\Sigma$ is embedded. 

To finish, we show that the elliptic data $(f,\Gamma)$ gives rise to an elliptic equation on the principal curvatures for $\Sigma$, in other words, $\Sigma$ satisfies the Maximum Principle. Therefore, the Alexandrov Reflection Method applied to $\Sigma$ will say that $\Sigma$ is part of a totally geodesic sphere, whose horospherical metric is given, up to a conformal diffeomorphism of $\mathbb S ^n$, by the standard metric $g_0$ on $\mathbb S ^n$.

\begin{remark}
{\it Note that, given any geodesic ball $B$ in $\mathbb{S}^n$, there exists a conformal diffeomorphism, $T$, of $\mathbb{S}^n$  such that $T(B)= \mathbb{S}^n_+$.}
\end{remark}

Nevertheless, we can go further, and we can deal with annular domains, as Escobar did \cite{E0} for the scalar curvature, in the fully nonlinear elliptic case. Let us denote by $\n \in \mathbb{S}^n_+ \subset \mathbb{S}^n$ the north pole and let $r < \pi/2$. Denote by $\mathbb{B}(\n , r)$ the geodesic ball in $\mathbb{S}^n$ centered at $\n $ of radius $r$. Note that, by the choice of $r$, $\partial \mathbb{S}^n_+ \cap \partial \mathbb{B}(\n ,r) = \emptyset$. 

Denote by $\mathbb A (r)= \mathbb{S}^n_+ \setminus \overline{\mathbb B (\n , r)}$ the annular region determined by $\mathbb{S}^n_+$ and $ \mathbb B (\n , r)$. Note that the mean curvature of $\partial  \mathbb B (\n , r) $ with respect to $g_0$ and the inward orientation along $\partial \mathbb A (r) $ is a constant $h(r)$ depending only on $r$. Therefore, our second task will be to find a conformal metric on $\mathbb A (r)$ satisfying an elliptic condition in the interior and whose boundary components $\partial  \mathbb B (\n , r)$ and $\partial \mathbb{S}^n_+$ are minimal.

In other words, given $(f, \Gamma)$ an elliptic data, find $\rho \in C^\infty (\mathbb A (r))$  
so that the metric $g=e^{2\rho}g_0 $ satisfies 
\begin{equation}\label{p2}
\left\{\begin{array}{rl}
f(\lambda(p))= 1,& \lambda(p)\in\Gamma, \, \, p\in \mathbb{A}(r),\\ \\

h(g)= 0, & \textrm{ on } \partial \mathbb{B}(\n ,r)\cup \partial \mathbb{S}^n_+ . 
\end{array} \right.
\end{equation}

In the above situation, we will obtain

\begin{theorem}\label{thp2}
Let $\rho \in C^\infty (\overline{\mathbb{A}(r)})$ be a solution to \eqref{p2}. Then, $g= e^{2\rho}g_0$ is rotationally symmetric metric on $\overline{\mathbb{A}(r)}$. 
\end{theorem}

The strategy here is as above. We will construct a compact embedded horospherically convex hypersurface $\Sigma$ which is a topological annulus $\mathbb S^{n-1} \times [0,1]$. As we did above, we will see that $\Sigma$ is contained in the slab determined by two totally geodesic hypersurfaces. Also, we will see that $\Sigma$ is orthogonal to the totally geodesic hypersurfaces at each component of the boundary. Hence, we can extend $\Sigma$ across the totally geodesic hyperplanes by reflections to obtain a properly embedded hyper surface with two points at infinity. Therefore, using the Alexandrov Reflection Method, we will conclude that $\Sigma$ is rotationally symmetric and so is $g$.

In \cite{CHY}, Chang-Han-Yang have classified all posible radial solution to $\sigma_k -$Yamabe problem, that is, to the equation $ \sigma _k (\lambda _1 , \ldots , \lambda _n) = 1.$
Therefore, when $f = \sigma _k$, we can conclude that the solution in Theorem \ref{thp2} is one given in \cite{CHY}. From the geometric point of view, the solution in Theorem \ref{thp2} corresponds either to Delaunay type hypersurfaces or totally umbillical spheres.

As Escobar pointed out \cite{E0,E1,E2} for the constant scalar curvature case, there are examples of rotational metrics with prescribed constant mean curvature on its boundary. {\it Does there exist a general classification result in this case?} Our technique for the fully nonlinear case uses strongly the assumption of minimality on its boundary, as well as Escobar's proof for the constant scalar curvature case. 

Up to this point, we have focused on dimension $n \geq 3$. Nevertheless, the two dimensional case is of special interest. Classically, the techniques in dimension $n=2$ are based on Complex Analysis, however, our method also works in this case. In fact, it is clear that the Schouten tensor is only defined for dimensions greater or equal than $3$, but, we can still define a symmetric two tensor for any conformal metric $g=e^{2\rho} g_0$ on the standard $2-$sphere $(\mathbb S ^2 , g_0)$. In particular, the trace of such symmetric tensor will be the Gaussian curvature of the conformal metric $g$ or, in other words, we will be dealing with solutions to the classical 
\emph{Neumann problem for the Liouville equation} in dimension 2. We will explain this in detail in the last Section and we will see in Theorems \ref{thp3} and \ref{thp4} how the previous results extend to this case.  So, these results will classify the space of solutions to Neumann problems for fully nonlinear conformally invariant equations, extending previous results for the classical Liouville equation (cf. \cite{GaMi,HaWa,J,Zha} and references therein). We shall remark the previous cited results rely strongly on the complex variable one can introduce in dimension $2$, we use the geometric approach developed in this paper. 

In a forthcoming paper \cite{AE}, the authors study the degenerate elliptic case.

\bigskip
\section{A Representation Formula}

In this section we describe the main tools that will be need to prove Theorems \ref{thp1} and \ref{thp2}.
Firstly, we will present the representation formula given in \cite{EGM} for domains of the round sphere 
endowed with a conformal metric as an immersion into the hyperbolic space.

Let $(\mathbb S ^n , g_0)$ be the standard $n-$sphere. Let $\overline \Omega\subset \mathbb S^n$  be a compact domain and $g=e^{2\rho}g_0$ be a $C^\infty$ metric on $\overline \Omega$. Assume that 
$$
{\rm Sch}_g(p)<\frac{1}{2} \text{ for all } p\in \overline \Omega,
$$
that is, each eigenvalue of the Schouten tensor is less than $1/2$. We can always achieve this 
condition by a dilation $g_t=e^{2t}g$, since $\overline \Omega$ is compact.

Denote by $\mathbb L ^{n+2}$ the standard Lorentz-Minkowski space, i.e, 
$\mathbb L ^{n+2} = (\real ^{n+2} , \meta{}{})$, where $\meta{}{}$ is the standard Lorentzian 
metric given by 
$$ \meta{}{} = - dx_0 ^2 + \sum _{i=1}^{n+1} dx_i^2 .$$

In this model one can consider
\begin{equation*}
\begin{array}{rcl}
\mathbb H ^{n+1} &=& \{  x \in \mathbb L ^{n+2} \, : \, \, \meta{x}{x}=-1 , \, x_0 >0\} ,\\
\mathbb S ^{n+1}_1 &=& \{  x \in \mathbb L ^{n+2} \, : \, \, \meta{x}{x}= 1\} , \\
\mathbb N ^{n+1}_+ &=& \{  x \in \mathbb L ^{n+2} \, : \, \, \meta{x}{x}= 0 , \, x_0 >0\} ,
\end{array}
\end{equation*}that is, the Hyperbolic Space, the deSitter Space and the Light Cone respectively.

Following \cite{BEQ}, one can construct a representation of $(\Omega, g)$ as an immersion 
$\phi:\overline \Omega \to \mathbb H^{n+1} \subset (\mathbb{L}^{n+2} , \meta{}{})$, endowed with a canonical orientation  $\eta : \overline{\Omega} \to \mathbb{S}^{n+1}_1 \subset \mathbb L ^{n+2}$, given by

\begin{equation}\label{phi}
\phi(x) = \frac{e^\rho}{2}\big(1+e^{-2\rho}(1+ \|\nabla^0 \rho  \|_0^2)\big)(1,x)+e^{-\rho}(0,-x+\nabla^0\rho)
\end{equation} and whose hyperbolic Gauss map is given by  $G(x)=x$.  In other words, one can construct a horospherically concave hypersurface $\Sigma=\phi(\Omega)$  with boundary $\partial\Sigma = \phi(\partial \Omega)$. 
Here, $\| \, \cdot \, \|_0$ and $ \nabla^0  \rho$ represent the norm and the gradient with respect to $g_0$. 

Recall that the hyperbolic Gauss map is defined as follows. Let $x \in \Omega$ be a point in our domain and consider $p := \phi (x) \in \mathbb H ^{n+1}$ and $v := -\eta (x) \in T_p \mathbb H ^{n+1}$. Then, $G : \Omega \to \mathbb S ^n $ is defined by 
$$ G( x ) := \lim _{t \to + \infty} \gamma _{p,v} (t)  \in \mathbb S ^n , $$where $\gamma : \mathbb R \to \mathbb H ^{n+1}$ is the complete geodesic parametrized by arc-length in $\mathbb H ^{n+1}$ passing through $p$ in the direction $v$. 

\begin{remark}
{\it Note that, from \eqref{phi}, the immersion is smooth. Moreover, one can see that $\phi $ is $C^k$, and also the First and Second Fundamental Forms are $C^{k-1}$, when $\rho$ is $C^{k+1}$. So, we could relax the differentiability hypothesis on $\rho $ in order that $\Sigma$ satisfies the Maximum Principle.} 
\end{remark}

\medskip 

We recall that  an immersion is horospherically concave if and only if the principal curvatures at any point are bigger than $-1$ for the prescribed orientation $\eta $ (the inward orientation for a totally umbilical sphere). Also $ g :=  \meta{d\psi}{d\psi} = e^{2\rho} g_0 $ is a Riemmanian metric, being $\psi := \phi - \eta : \Omega \to \mathbb N ^{n+1}_+$ and satisfies 
$$\psi = e^{\rho}(1,x), \, x\in \overline{\Omega},$$that is, $g$ is nothing but the First Fundamental Form of $\psi$. Moreover, the principal curvatures, $\kappa _i $, of 
$\Sigma $ and the eigenvalues, $\lambda _i$, of the Schouten tensor of $g$ are related by
\begin{equation}\label{lambdakappa}
\lambda _ i = \frac{1}{2} - \frac{1}{1+ \kappa _i} .
\end{equation}

\bigskip

\section{Proof of Theorem \ref{thp1}}\label{s3}

As we discussed in the Introduction, we will divide the proof in parts. Assume $g = e^{2\rho } g_0$ is a solution 
of problem \eqref{p1}. 


\vspace{.3cm}

\begin{quote}
{\bf Claim A: }{\it 
There exists $t_0 >0 $ so that the dilated metric $g _{t_0} := e^{2t_0} g$ satisfies: 
\begin{itemize}
\item $g_{t_0} $ is  a solution of the elliptic problem
\begin{equation*}
\left\{\begin{array}{rl}
f_{t_0}(\lambda ^{t_0}(p))= 1,& \lambda ^ {t_0}(p)\in\Gamma _{t_0}, \, \, p\in \mathbb{S}^n_+,\\ \\

h(g_{t_0}) = e^{-t_0}c	, & \textrm{ on } \partial \mathbb{S}^n_+ .
\end{array} \right. 
\end{equation*}where $f_{t_0}(\lambda ^{t_0}(p)) = f(e^{-t_0}\lambda (p))$ and $\Gamma _{t_0} = e^{-t_0}\Gamma$. \\[3mm]

\item $| {\rm Sch}_{g_{t_0}} | < 1/2$. \\[3mm]
\item $ \dfrac{1+2\lambda _i ^{t_0}}{1-2\lambda _i ^{t_0}} > |\tanh ( \sinh ^{-1}(h(g_{t_0})))|$. 
\end{itemize}}
\end{quote}
\begin{proof}[Proof of Claim A]
Consider a dilation of $g$, i.e. $g_t=e^{2t}g$. Hence, the eigenvalues of the Schouten tensor of $g_t$, denoted by $\lambda _i ^t$, change by the formula
$$ \lambda _i ^t = e^{-t} \lambda _i ,$$where $\lambda _i $ are the eigenvalues of the Schouten tensor of $g$. 

Moreover, when we dilate the metric $g$ the mean curvature of the boundary changes as
$$
h_t = h(g_t) = e^{-t}c ,
$$this proves the first item.

Now, since $\overline{\mathbb{S}^n_+}$ is compact and $\lambda _i ^t = e^{-t}\lambda _i$, there exists $t_1 >0 $ such that $ |\lambda _i ^t | < 1/2$ for all $i =1, \ldots ,n$ and $t \geq t_1$.

Since
$$ \frac{1+2 \lambda^t _i}{1-2\lambda^t _i} =\frac{1+2 e^{-t}\lambda _i}{1-2e^{-t}\lambda^t _i} \to 1 \text{ as } t \to +\infty$$and 
$$ h(g_{t})=e^{-t} c \to 0 \text{ as } t \to +\infty ,$$there exists $t_2 >0$ such that 

$$ \frac{1+2 \lambda^t _i}{1-2\lambda^t _i} > |\tanh(\sinh  ^{-1} (h(g_t)))| $$for all $t \geq t_2$. 

Thus, taking $t _0 \geq {\rm max}\{t_1 ,t_2 \}$, we prove the second and third items.
\end{proof}

Hence, from this point on, we can assume that $g$ is a solution of \eqref{p1} and satisfies
\begin{itemize}
\item[(P1)] $| {\rm Sch}_{g} | < 1/2$. \\[1mm]
\item[(P2)] $ \dfrac{1+ 2\lambda _i }{1-2\lambda _i } > |\tanh ( \sinh  ^{-1}(h(g)))|$ for all $i = 1, \ldots ,n$.
\end{itemize}

Let us explain geometrically the meaning of the above properties. (P1) will say that we can construct an immersion $\phi : \Omega \to \mathbb H ^{n+1}$ from the conformal factor $\rho$ (cf. \cite{EGM}). As we saw above, the principal curvatures of $\phi$ are related to the eigenvalues of the Schouten tensor as 
$$ \kappa _i =   \dfrac{1+ 2\lambda _i }{1-2\lambda _i }, \,\text{ for all }� i = 1, \ldots , n , $$and hence, (P2) implies that 
$$ \kappa _i > |\tanh ( \sinh  ^{-1}(h(g)))|, \text{ for all } i = 1, \ldots ,n .$$

So, we can rewrite 
\begin{itemize}
\item[(P2')] $ \kappa _i> |\tanh ( \sinh  ^{-1}(h(g)))|$ for all $i = 1, \ldots ,n$.
\end{itemize}

The above properties for $g$ will lead us to the following geometric construction associated to $(\overline{\mathbb{S}^n_+}, g)$:

\vspace{.3cm}

\begin{quote}
{\bf Claim B: }{\it Given a conformal metric on $\overline{\mathbb S ^n _+}$, there exists a horospherically concave (in fact convex) hypersurface $\Sigma \subset \mathbb H ^{n+1}$ with compact boundary $\partial \Sigma$ such that $\Sigma $ and $\partial \Sigma$ are topologically $\mathbb S ^n _+$ and $\partial \mathbb S^n_+ = \mathbb S ^{n-1}$ respectively.

Moreover, $\partial \Sigma$ is contained on an equidistant hypersurface at distance $|\sinh ^{-1}(c)|$ 
(the orientation depending on the sign of $c$) to a totally geodesic hyperplane.}
\end{quote}
\begin{proof}[Proof of Claim B]
Now, consider $g=e^{2\rho}g_0$ in $\overline {\mathbb{S}^n_+}$ satisfying Property (P1) above.
From \cite{BEQ} we have that $\phi$ given by (\ref{phi}) is a smooth injective 
immersion of $\mathbb S^n_+$ onto 
$\mathbb H ^{n+1}$  and by applying the above construction and get $\Sigma = \phi (\mathbb S ^n _+)$ and $\partial \Sigma = \phi (\partial \mathbb S ^n _+)$.  This proves the first part of the Claim.

Now, we will study the behavior of the boundary. Since $\partial \mathbb{S}^n_+ $ has constant mean curvature $c$ with respect to $g$, we have that $\rho$ satisfies the boundary equation 
\begin{equation*}\label{}
e^{-\rho} \frac{\partial \rho}{\partial \nu}= c\, \,\textrm{ on } \partial\mathbb{S}^n_+ .
\end{equation*}

Considering the canonical coordinates 
$$(x_1,\ldots,x_{n+1})\in\mathbb S^n_+\subset \real^{n+1} = \{ (x_0, x_1, \ldots , x_{n+1})
 \in \real ^{n+2} \, : \, \, x_0 =1\} \subset \mathbb L ^{n+2},$$
we have that the inward unit normal vector field of $\mathbb S^n_+$ is just $\nu=e_{n+1}$. 
We note that the above condition on the mean curvature says that the vector
$$
Y:= e^{-\rho}\left( -x+\nabla^0\rho \right) = e^{-\rho}(-x_1+\frac{\partial \rho}{\partial e_1},
\ldots,-x_{n+1}+\frac{\partial \rho}{\partial e_{n+1}})
$$
satisfies
$$
Y_{|_{\partial \mathbb{S}^n_+}} = 
\big(e^{-\rho}(-x_1+\frac{\partial \rho}{\partial e_1}), \ldots,e^{-\rho}(-x_{n}+\frac{\partial \rho}{\partial e_{n}}),c\big), 
$$
that is, $ Y(x)\in \{(x_1,\ldots,x_n,x_{n+1})\in\real^{n+1} \, : \, \, x_{n+1} =c\}$ for all  $x\in\partial   \mathbb{S}^n_+ $, since $x_{n+1} =0 $ along $\partial   \mathbb{S}^n_+$.

Set $P(c) := \{ (x_0, x_1 , \ldots , x_{n+1}) \in \real ^{n+2}  \, : \, \, x_{n+1} =c \}$. 
Thus, from \eqref{phi} and the above observation we see that  
$$  
\phi (x)  \in \mathbb H ^{n+1} \cap P(c) \text{ for all } x \in \partial \mathbb S ^n _+.
$$
In other words, $\partial \Sigma=\phi (\partial  \mathbb{S}^n_+ )\subset E(c)$ where $ E(c) := \mathbb H^{n+1}\cap P(c)$ is the equidistant hypersurface at distance $|\sinh ^{-1}(c)|$ to a totally geodesic hyperplane.
In particular $\partial \Sigma$ lies in the totally geodesic hyperplane when $c=0$, 
i.e., a minimal boundary for the conformal metric.
\end{proof}

So, let us explain the above Claim B in a different model of the Hyperbolic Space. Consider the Poincar\'{e} Ball model $(\mathbb B ^{n+1} , g_{-1})$ of $\mathbb H ^{n+1}$, where $\mathbb B ^{n+1}$ is the Euclidean ball in $\real ^{n+1}$ of radius $1$ and $g_{-1}$ is the Poincar\`{e} metric. It is well known that the boundary at infinity of $\mathbb H ^{n+1}$ in the Poincar\'{e} Ball model corresponds to 
$$ \partial _\infty \mathbb H ^{n+1} = \partial \mathbb B ^{n+1} = \mathbb S ^n ,$$and we can consider $\mathbb S ^n _+$ as the upper hemisphere here and $\partial \mathbb S ^n _+$ as the equator 
$$ \partial  \mathbb S ^n _+ = \{ (x_1, \ldots , x_{n+1}) \in \mathbb S^{n} \, : \, \, x_{n+1}=0\} \subset \real ^{n+1} .$$

The totally geodesic hyperplane with boundary at infinity $\partial \mathbb S ^n_+$ is given by 
$$ E(0) := \{ (x_1 , \ldots , x_{n+1}) \in \mathbb B ^{n+1} \, : \, \, x_{n+1} =0 \}  $$and the equidistant at distance $ \sinh^{-1} (c) $ to $E(0)$ is given by 
$$ E(c)  = \{ {\rm exp }_p (\sinh ^{-1}(c) \, N (p)) \in \mathbb B ^{n+1}  \, : \, \, p \in E(0)\}  ,$$where ${\rm exp}$ is the exponential map in $\mathbb H ^{n+1}$ and $N$ is the upward normal along $E(0)$, i.e., the one pointing in the half-space containing containing $\mathbb S^n_+$ as its boundary at infinity. Recall that each $E(c)$ is totally umbilic with constant principal curvatures $-\tanh (\sinh ^{-1}(c))$, that is, 
$$ II_{E(c)} = -\tanh (\sinh ^{-1} (c)) I _{E(c)},$$where $I_{E(c)}$ and $II_{E(c)}$ denote the First and Second Fundamental Form respectively. Here the orientation for $E(c)$ the one pointing at the component of $\mathbb H ^{n+1} \setminus E(c)$ containing the north pole ${\bf n}$ in its boundary at infinity.

Without loss of generality we will assume $c \geq 0$.

\vspace{.3cm}

\begin{quote} 
{\bf Claim C: }{\it $\Sigma \subset E(c)^+$, where $E(c)^+$ is the half-space containing $\mathbb S^n_+$ as its boundary at infinity.} 
\end{quote}
\begin{proof}[Proof of Claim C]
Assume we have points $p\in \Sigma \cap E(c)^-$. Since $\overline\Sigma$ is compact and $\partial \Sigma \subset E(c)$ we can find a lowest point $p_0$. 

First, note that $E(s)$, $ s \in (-\infty ,c)$ is a foliation (by totally umbilical hypersurfaces) of $E^- (c)$. For $s$ close to $-\infty$, $E(s)$ is disjoint from $\overline \Sigma$. So, we start increase $s$ up to $c$. Suppose we find the first contact point  of $\overline \Sigma$ with some $E(\bar s)$, $\bar s\in (-\infty, c)$. 

Then, at this first contact point $p_0\in E(\bar s)\cap \Sigma$, the canonical orientation $\eta$ of $\Sigma$ at $p_0$ and the unit normal of $E(\bar s)$ pointing upwards coincide, this follows from property (P2). 
More precisely, condition (P2') says that the principal curvatures $\kappa _i$, computed with respect to the canonical orientation, satisfy
\begin{equation}\label{Eq:Relation}
\kappa _i > k _i (s) \text{ for all }  s \in (-\infty ,c] \text{ and } i = 1, \ldots , n, 
\end{equation}where $k_i (s)$, $i=1, \ldots ,n$, are the principal curvatures of $E(s)$ computed with respect to the orientation  pointing at the component of $\mathbb H ^{n+1} \setminus E(s)$ containing the north pole ${\bf n}$ in its boundary at infinity. So, at the first contact point, $\Sigma $ is locally at one side of $E(\bar s)$, $\Sigma$ must be locally contained on $E(\bar s)^+$. Since $\Sigma$ is more curved than $E(\bar s)$ at $p_0 \in \Sigma$, we obtain that the canonical orientation must coincide with the upwards orientation. This implies that $G(p_0)\in \mathbb S^n_-$, which contradicts $G(\Sigma)=\mathbb{S}^n_+$.

Moreover, the only possibility is $\overline{\Sigma} \cap E(c) = \partial \Sigma$. In fact, if there were another intersection point $ p_0 \in \overline{\Sigma} \subset E(c) \setminus \partial \Sigma$, the same argument above apply and hence we will get a contradiction.

\end{proof}

Next, we will study how $\Sigma$ intersects $E(c)$.

\vspace{.3cm}

\begin{quote} 
{\bf Claim D: }{\it $\Sigma$  makes a constant angle 
$$ \alpha(c) := {\rm arccos}\left(  \frac{c}{\sqrt{1+c^2}}\right) $$ with $E(c)$ along $\partial \Sigma$. Here 
$ \alpha (c)  $ is the angle between the canonical orientation $\eta$ and the upward normal along $E(c)$.}
\end{quote}
\begin{proof}[Proof of Claim D]
From the definition of the immersion $\phi$  one easily checks that its unit normal is given by
$$
\eta(x) = \frac{e^{-\rho}}{2} \big(\| \nabla \rho\|^2 -1-e^{2\rho}  \big) (1,x) + e^{-\rho}(0,-x+ \nabla \rho).
$$

Therefore

$$
\eta_{|_{\partial \mathbb S^n_+}} \in \mathbb S^{n+1}_1\cap \{(x_0,\ldots, x_n, c)\in\mathbb L^{n+2} \} = \mathbb S^{n+1}_1\cap P_c.
$$

The normal along $E(c)= \mathbb H ^{n+1} \cap P(c)$ is given by 
$$ n(x) = \frac{1}{\sqrt{1+c^2}} e_{n+1} + \frac{c}{\sqrt{1+c^2}} x, \, \, x \in E(c)= \mathbb H ^{n+1} \cap P_c,$$and so
\begin{equation}\label{angle}
\meta{\eta (x)}{n(\phi (x))} = \frac{c}{\sqrt{1+c^2}}\, \text{  along } \partial \Sigma . 
\end{equation}

Therefore, $\Sigma$ makes a constant angle $ {\rm arccos}\left(  \frac{c}{\sqrt{1+c^2}}\right) $ with $E(c)$ along $\partial \Sigma$.

\end{proof}

\medskip 

\begin{quote} 
{\bf Claim E: }{\it $\Sigma $ is embedded.}
\end{quote}
\begin{proof}[Proof of Claim E]
Note that $E(c)$ is isometric to a Hyperbolic space with constant sectional curvatures $K_{sect} (t) = \tanh ^2 (\sinh ^{-1}(c)) -1$. Moreover, as we already saw, (P2') says that the principal curvatures $\kappa _i$, computed with respect to the canonical orientation, satisfy
\begin{equation*}
\kappa _i > k _i (c) \text{ for all }   i = 1, \ldots , n, 
\end{equation*}where $k_i (c)$, $i=1, \ldots ,n$, are the principal curvatures of $E(c)$ computed with respect to the orientation  pointing at the component of $\mathbb H ^{n+1} \setminus E(c)$ containing the north pole ${\bf n}$ in its boundary at infinity. From this, we conclude that $\partial \Sigma \subset E(c)$ is convex and compact in $E(c)$. 
%
In fact, let $\alpha : (-\epsilon , \epsilon) \to \partial \Sigma \subset E(c)$ be a curve along the boundary. We must compute the second fundamental form of $\partial \Sigma$ in $E(c)$, that is 
$$ II_{\partial \Sigma} (\alpha ' ,\alpha ')= \meta{\nabla ^{E(c)} _{\alpha '} \alpha '}{N} ,$$where $\nabla ^{E(c)}$ is the connection on $E(c)$ and $N$ is the inward normal along $\partial \Sigma$ in $E(c)$. 

\begin{remark}
{\it Note that here we use the geometric definition for the second fundamental form, opposite to the usual analytic definition used at the begging when we explain the Yamabe Problem with boundary. They differ by a sign.}
\end{remark}

From \eqref{angle}, one can compute $N$ as 
$$ N := \sqrt{1+c^2} \, \eta - c \, n ,$$where $\eta $ is the normal along $\Sigma$ and $n$ is the normal along $E(c)$. Therefore, from (P2'), we obtain
\begin{equation*}
\begin{split}
II_{\partial \Sigma}(\alpha ',\alpha ') & = \meta{\nabla ^{E(c)} _{\alpha '} \alpha '}{N}= \meta{\nabla ^{\mathbb H^{n+1}} _{\alpha '} \alpha '}{N} = \\
 &= \sqrt{1+c^2} \meta{\nabla ^{\mathbb H^{n+1}} _{\alpha '} \alpha '}{\eta} - c \meta{\nabla ^{\mathbb H^{n+1}} _{\alpha '} \alpha '}{n} \\
 & =\sqrt{1+c^2} \, II_{\Sigma}(\alpha ',\alpha ') - c \, II_{E(c)}(\alpha ',\alpha ') > 0,
\end{split}
\end{equation*}which proves the claim. $\partial \Sigma $ is convex in $E(c)$.
%
%
Hence, by Do Carmo-Warner Theorem \cite{CW}, $\partial \Sigma$ is embedded in $E(c)$.

Consider the foliation of $E(c)^+$ given by the translations $E^s(c)$ of $E(c)$, with $s\geq 0$ and $E^0(c) = E(c)$. We claim that when the intersection $\Sigma_s=\Sigma \cap E^s(c)$ is transversal, each connected component of $\Sigma \cap E^s (c)$ is a compact and convex hypersurface in $E^s (c)$ and therefore embedded (see \cite{CW}). 
In fact, arguing as above, we consider a smooth curve
%
 $\alpha : (-\epsilon , \epsilon) \to \Sigma \cap E^s(c)$ parametrized by arc-length. Then
$$ \nabla ^s _{\alpha '} \alpha ' + II_s(\alpha ' , \alpha ') n_s= \nabla ^{\mathbb H ^{n+1}}_{\alpha '} \alpha ' = \nabla ^\Sigma _{\alpha '} \alpha ' + II_\Sigma(\alpha ' , \alpha ')\eta.$$

Thus,
$$\| \nabla ^s _{\alpha ' } \alpha ' \|^2 = \|  \nabla ^\Sigma _{\alpha ' } \alpha ' \|^2 + II_\Sigma (\alpha '  , \alpha ' ) ^2 - II_s (\alpha '  , \alpha ' )^2 >0, $$so, each transverse intersection is convex in $E^s (c)$.

For $s $ close to $0$, say $0\leq s\leq \varepsilon$, the intersection $\Sigma \cap E^s(c)$ is transverse and $\Sigma _s$ has only one component
which is embedded and homeomorphic to $\mathbb S^{n-1}$.  
Denoting by $\Sigma_0$ the portion of $\Sigma$ below  $E^\varepsilon(c)$ we have that $\Sigma_0$ is an embedded annulus.
Now we can glue $\Sigma\setminus  \Sigma_0$ with a embedded convex ball in order to construct a closed convex hypersurface (at least $C^2$) $\Sigma_1$ in $\mathbb H^{n+1}$. Hence, by Do Carmo-Warner Theorem, we conclude that $\Sigma_1$ is embedded and thus $\Sigma$ is embedded.

We also can prove that $\Sigma$ is embedded showing that $\Sigma _s$ has only one component or it is empty for all $s >0 $ by using the condition \eqref{Eq:Relation} and Do Carmo-Warner Theorem.

\end{proof}

So, up to a dilation on $g$, we can assume that $\Sigma$ is embedded,  convex with respect to the inward 
orientation, compact (homeomorphic to $\overline {\mathbb S^n_+}$),  with boundary $\partial \Sigma \subset E(c)$ that makes a constant angle with $E(c)$, and $\Sigma \subset  E(c)^+$.

\vspace{.3cm}

\begin{quote}
{\bf Claim F: }{\it $\Sigma$ satisfies an elliptic equation $(\mathcal W , \Gamma ^*)$ for a given curvature function depending on its principal curvatures; that is, 

\begin{equation}\label{p1W}
\left\{\begin{array}{rl}
\mathcal W (\kappa(p))= 1,& \kappa (p)\in\Gamma ^*, \, \, p\in \Sigma,\\ \\

\angle( \partial \Sigma  , n) = \alpha (c) &\text{ along } \,  \partial \Sigma \subset E(c) ,
\end{array} \right. 
\end{equation}where, 
$$ \angle( \partial \Sigma  , n) = \alpha (c)$$means that $\partial \Sigma \subset E(c)$ makes a constant angle $\alpha (c)= {\rm arccos}\left(\frac{c}{\sqrt{1+c^2}} \right)$ with $E(c)$. Here, $\kappa (p)$ denotes the principal curvature vector, i.e., $\kappa (p)= (\kappa _1 (p), \ldots , \kappa _n (p))$.}
\end{quote}
\begin{proof}[Proof of Claim F]
We recall from \cite[Section 4]{BEQ} the definition of elliptic data for a horospherically concave hypersurface in $\mathbb H ^{n+1}$.  Let 
$$\Gamma ^* _n =\{ (x_1, \ldots , x_n) \in \real ^n \, : \,\, x_i >1 \}$$and
$$\Gamma ^* _1 =\{ (x_1, \ldots , x_n) \in \real ^n \, : \,\, \sum _{i=1}^nx_i >n  \}.$$

Consider a symmetric function $\mathcal W (x_1 , \ldots , x_n)$ with $\mathcal W (1, \ldots , 1) = 0$ and $\Gamma ^*$ an open connected component of 
$$ \{ (x_1, \ldots , x_n) \in \real ^n \, : \,\, \mathcal W (x_1, \ldots ,x_n ) >0  \}. $$

We say that $(\mathcal W , \Gamma ^*) $ is an elliptic data if they satisfy

\begin{enumerate}

\item[(1)] $\Gamma ^*_n \subset \Gamma ^* \subset \Gamma ^* _1$; \\

\item[(2)] For all $(x_1, \ldots , x_n)\in \Gamma ^*$ and all $(y_1, \ldots , y_n) \in \Gamma ^* \cap ((x_1,\ldots ,x_n) + \Gamma^* _n)$, there exists a curve $\gamma$ connecting $(x_1, \ldots , x_n)$ to $(y_1, \ldots , y_n)$ inside $\Gamma^*$ such that $\gamma ' \in \Gamma ^*_n $ along $\gamma$; \\

\item[(3)] $\dfrac{\partial \mathcal W}{\partial x_i}>0$, \quad $\forall \, i=1\ldots, n.$

\item[(4)] There exists $r_0>0$ such that $\mathcal W (r_0, \ldots , r_0) =1$.
\end{enumerate}

The ellipticity of $(f, \Gamma)$ implies that, if we define 
$$ \mathcal W (\kappa (p)) = (f \circ T )( \kappa (p)) ,$$where
$$ T(x_1, \ldots, x_n) = \left( \frac{x_1-1}{2(x_1 +1)} ,\ldots, \frac{x_n-1}{2(x_n +1)} \right),$$then $(\mathcal W, \Gamma ^*)$, where $\Gamma ^* = T^{-1}(\Gamma)$, is elliptic.
\end{proof}

We shall do a couple of remarks here. The ellipticity means that the horospherically concave hypersurface we construct from the conformal metric $g= e^{2\rho} g_0$ satisfies the Maximum Principle (Interior and Boundary). Moreover, the above condition (4) in the ellipticity of $(\mathcal W, \Gamma ^*)$ means that there exists a totally umbilical sphere $S_0$ so that its principal curvatures satisfy the equation $\mathcal W (\kappa (p)) =1$ for all $p \in S_0$.

\vspace{.3cm}

\begin{quote}
{\bf Claim G: }{\it $\Sigma$ is rotationally symmetric.}
\end{quote}
\begin{proof}
Up to this point, $\Sigma $ is embedded, convex , compact (homeomorphic to $\overline {\mathbb S^n_+}$),  with boundary $\partial \Sigma \subset E(c)$ that makes a contant angle $\alpha (c)$ with $E(c)$, and $\Sigma \subset  E(c)^+$. Moreover, Claim F implies that $\Sigma $ satisfies the (Interior and Boundary) Maximum Principle (see \cite[Section 4]{BEQ} and \cite{K}). Let $D \subset E(c)$ be the domain bounded by $\partial \Sigma $ in $E(c)$ and set $\Omega \subset \mathbb H ^{n+1}$ the domain bounded by $\Sigma \cup D$. 


So, Claim G follows from an application of the Alexandrov Reflection Principle. Take any geodesic $\beta \subset E(0)$, parametrized by arc-length $s$, and let $P(s)$ the foliation by totally geodesic hyperplanes orthogonal to $\beta$ at $\beta (s)$, let $P^+ (s)$ (resp. $P^- (s)$) denote the half-space determine by $P(s)$ such that $\lim _{s\to +\infty } \beta (s) \in \partial _\infty P^+ (s)$ (resp. $\lim _{s\to -\infty } \beta (s) \in \partial _\infty P^- (s)$). Moreover, let us denote by $R_s$ the reflection (which is an isometry) throughout $P(s)$. Such isometry $R_s$ leaves invariant $P(s)$ pointwise.

For $s$ close to $-\infty$, $P(s) \cap \Sigma =\emptyset $. So, we increases $\bar s$ up to the first contact point of $P(\bar s)$ and $\Sigma $. For $s >\bar s$ , close enough to $\bar s$, the part $\Sigma ^- (s) = \Sigma \cap P^- (s)$ is a graph (with boundary) over $P(s)$ and $\tilde \Sigma ^+ (s)= R_ s (\Sigma ^- (s)) \subset \Omega$. One important point to note here is that $R_s(E(h_t)) = E(h_t)$ for all $s \in \mathbb R$, that is, $E(h_t)$ is invariant by $R_s$ for all $s \in \mathbb R$. This implies that $R_s(\partial \Sigma) \subset E(h_t)$, i.e., the reflection $R_s(p)$ of any point $p \in \partial \Sigma \subset E(h_t)$ can not belong to the interior of $\Sigma$, that is, 
\begin{equation}\label{Eq:Contact}
R_s(p)\not\in \Sigma \setminus \partial \Sigma \text{ for all } p \in \partial \Sigma \subset E(h_t) .
\end{equation}

So, we can increase $s$ up to the first contact point of $\tilde \Sigma ^+ (s)$ and $\Sigma ^+ (s) = \Sigma \cap P^+(s)$, such a first contact point exists since $\Sigma$ is compact. Moreover, such first contact point must be either 
interior or boundary point by \eqref{Eq:Contact}. So, in any case, by the Maximum Principle (Interior or Boundary), we obtain $\tilde \Sigma ^+ (\bar s)= \Sigma ^+ (\bar s)$ for some $\bar s \in \mathbb R $, which means that $P(\bar s)$ is a plane of symmetry for $\Sigma$. Observe that for applying the Maximum Principle if the first contact point occurs at the boundary, we have strongly used the fact that $\Sigma$ makes a constant angle with $E(c)$ along $\partial \Sigma$.

Thus, we can perform the Alexandrov Reflection Method for any geodesic contained in $E(0)$ and hence $\Sigma$ must be rotationally symmetric. 

\end{proof}

Let $\gamma \subset \mathbb H ^{n+1}$ be the complete geodesic joining the south ${\bf s}$ and the north pole ${\bf n}$. Since $\Sigma $ is rotationally symmetric, up to an isometry, we can assume that $\gamma$ is the axis of rotation. We must point out that the isometry $i : \mathbb H ^{n+1} \to \mathbb H ^{n+1}$ that we use is just a hyperbolic translation whose fixed points belong to $\partial \mathbb S ^n _+$. Moreover, it is well-known that any isometry of $\mathbb H ^{n+1}$ induces an unique conformal diffeomorphism of $\mathbb S ^n$ and viceversa. So, note that the above hyperbolic translation induces a conformal diffeomorphism that leaves invariant $\mathbb S ^n _+$. Therefore, the image of the hyperbolic Gauss map of $\tilde \Sigma = i (\Sigma)$  remains invariant, i.e., $\tilde G (\tilde \Sigma) = G (\Sigma) = \mathbb S ^n _+$. 

From \cite[Lemma 3.2; Section 3]{E}, we know that if the associated horospherically concave hypersurface to a conformal metric  is invariant under a subgroup of isometries of $\mathbb H ^{n+1}$, the conformal metric is invariant under the subgroup of conformal diffeomorphism induced by the subgroup of isometries, and viceversa. Therefore, $g$ is rotationally symmetric. 

Since the problem we are considering is conformally invariant we can summarize what we have done as:

\vspace{.2cm}

\begin{quote}
{\it Given $\rho \in C^\infty (\overline{\mathbb S ^n _+})$ a solution to \eqref{p1}, up to a dilation and a conformal diffeomorphism that leaves invariant $\mathbb S ^n _+ $, the conformal metric $g = e^{2\rho } g_0$ induces an embedded rotationally symmetric hypersurface in $\mathbb H ^{n+1}$, horospherically concave, compact (homeomorphic to $\overline {\mathbb S^n_+}$),  with boundary $\partial \Sigma \subset E(c)$ that makes a constant angle $\alpha (c)$ with $E(c)$, and $\Sigma \subset  E(c)^+$. Therefore, $g$ is rotationally symmetric.}
\end{quote}

\vspace{.2cm}

Thus, at this point, we have shown that $g$ is rotationally symmetric. In what follows let us see that $g$ is in fact, 
up to a conformal diffeomorphism, the standard metric $g_0$. 

As we remarked above, condition (4) means that there exists a totally umbilical sphere $S_0 \subset \mathbb H ^{n+1}$ such that $\mathcal W (\kappa (p)) = 1$ for all $p \in S_0$.

\vspace{.3cm}

\begin{quote}
{\bf Claim H: }{\it $\Sigma $ is part of a totally umbilical sphere $S_0$}
\end{quote}
\begin{proof}
Let $\gamma \subset \mathbb H ^{n+1}$ be the complete geodesic joining the south ${\bf s}$ and the north pole ${\bf n}$. By the above discussion, we can assume $\gamma$ is the axis of symmetry for $\Sigma$. 

Let us denote by $S_0 (s)$ the totally umbilical sphere $S_0$ whose center is $\gamma (s)$. Since $S_0$ is compact and symmetric, it is clear that there exist $s_1 < s_2$ such that 
\begin{itemize}
\item $S_0 (s) \cap E(c) = \emptyset $  for all $s \in \mathbb R \setminus [s_1 ,s_2]$,
\item $S_0 (s) \cap E(c) $ is compact and transversal, in fact, it is a totally umbilical $(n-1)-$sphere in $E(c)$, for all $s \in ( s_1 ,s_2)$,
\item $S_0 (s_i) \cap E(c) = \{q\}$, i.e., the intersection is tangential. Note that $q\in \gamma$, $i=1,2$.
\end{itemize}

Let $s \in (s_1, s_2)$ and denote by $\nu (s)$ the inward normal vector field of $S_0 (s) \cap E(c) $ in $\mathbb{H}^{n+1}$, and $n$ normal along $ E(c) $ pointing upwards. Set 
$$ \Theta (s) = \meta{n}{\nu (s)} ,$$and note that $\Theta$ only depends on $s$ by the rotational symmetry of the problem. 

So, one can observe that 
\begin{itemize}
\item $\Theta (s) \to -1$ as $ s \to s_1$ and $s_1 < s$,
\item $\Theta (s) \to +1$ as $ s \to s_2$ and $s < s _2$,
\end{itemize}therefore, since one can observe that $\Theta$ is strictly decreasing, there exists $\bar s \in (s_1, s_2)$ such that $\Theta (\bar s) =\frac{c}{\sqrt{1+c^2}}$, that is, $S_0 (\bar s)$ makes a constant angle $\alpha (c)$ with $E(c)$. Set $S_0 ^+ (\bar s) = S_0 (\bar s) \cap E^+(c)$.

Since $\Sigma $ and $S_0 ^+ (\bar s)$ are rotationally symmetric it is easy to conclude that $\partial \Sigma = \partial S_0 ^+ (\bar s)$, otherwise we will get a contradiction with the Maximum Principle. Moreover, the problem reduces to a second order ODE by rotational symmetry, and both $\Sigma$ and $S_0 ^+ (\bar s)$ has the same initial conditions, i.e., the boundary is the same and makes the same angle with $E(c)$, we conclude that 
$$ \Sigma = S_0 ^+ (\bar s) ,$$as desired.
\end{proof}

So, it is easy to conclude that the horospherical metric of a totally umbilical sphere is (up to a conformal diffeormorphism) the standard metric $g_0$ on $\mathbb S ^n$ (see \cite{BEQ}). In fact, if the center of the totally umbilical is the origin (in the Poincar\'{e} Ball Model) then the horospherical metric is (up to possibly a dilation) the standard metric $g_0$. The dilation depends on the radius of the sphere, if we normalize so that $f(\lambda _1 , \ldots , \lambda _n) = 1$, or equivalently $\mathcal W (\kappa _1, \ldots , \kappa _n) = 1$, we will have an unique radius so that $ \mathcal W (\kappa (p))  =1$ for all points in this totally umbilical sphere. 

Therefore, since the center of $S_0 ^+ (\bar s)$ is $\gamma (\bar s)$, this means that $S_0 ^+ (\bar s) = T _{\bar s} (S_0 (0))$, where $S_0 (0)$ is the totally umbilical sphere whose center is the origin, and $T_{\bar s}$ is the hyperbolic translation along $\gamma $ at distance $\bar s$. Therefore, $g$ is the horospherical metric of $S_0 ^+ (\bar s)$ restricted to $\mathbb S ^n _+$, that is, $g=\Phi ^* g_0$ on $\mathbb S ^n _+$, where $\Phi$ is a conformal diffeomorphism of $(\mathbb S ^n ,g_0)$ associated to $T_{\bar s}$ as claimed.

\section{Proof of Theorem \ref{thp2}}\label{s4}

The proof of Theorem \ref{thp2} will follow the lines of Theorem \ref{thp1} but we must clarify certain aspects in this case. First, we shall introduce more notation. 

Let us denote by 
$$ \mathbb S (r) := \partial \mathbb B ({\bf n}, r) \subset \mathbb S ^n$$the boundary of the geodesic ball $\mathbb B ({\bf n}, r) $ in $(\mathbb S ^n , g_0)$ centered at the north pole ${\bf n}$ of radius $ r \in (0, \pi /2]$. Note that $r = \pi /2$ corresponds to the case $\mathbb B ({\bf n}, \pi /2) = \mathbb S ^n _+ $ and $\partial \mathbb B ({\bf n}, \pi /2 ) = \partial \mathbb S ^n _+  $.


For each $r \in (0 , \pi/2]$ there exists an unique totally geodesic hyperplane $E(r) \subset \mathbb H ^{n+1}$ whose boundary at infinity is $\partial _\infty E(r) = \mathbb S (r)$. When $r =\pi /2$, we denote $E= E(\pi/2)$, in this case, $\partial _\infty E = \mathbb S ^n _+$.

As we did above, and now that we know how dilations on $g$ works on the construction of the associated hypersurface $\Sigma$, we will do a dilation at the beginning and not to worry more about such dilation. 

Let $g$ be a conformal metric to $g_0$ defined on $\overline{\mathbb A (r)}$ satisfying \eqref{p2}. Choose $t \in \mathbb R $ such that
\begin{equation*}
\left| \frac{1+ 2e^{-t}\lambda _i (p)}{1-2e^{-t} \lambda _i (p)}\right|  > 0 \text{ for all } p \in \overline{\mathbb A (r)} \text{ and } i =1, \ldots , n.
\end{equation*}

So, from now on, we will work with the metric $g_t = e^{2t} g_0$, and we still denote it by $g$. Therefore, up to a dilation, we can assume that the eigenvalues of the Shouten tensor of $g$ satisfies
\begin{equation}\label{Eq:PrincipalAnnulus}
\left| \frac{1+ 2\lambda _i (p)}{1-2\lambda _i (p)}\right|  > 0  \text{ for all } p \in \overline{\mathbb A (r)} \text{ and } i =1, \ldots , n.
\end{equation}

\vspace{.3cm}

\begin{quote}
{\bf Claim A: }{\it Given $g$ as above, there exists a convex hypersurface $\Sigma \subset \mathbb H ^{n+1}$ with compact boundary $\partial \Sigma$ such that $\Sigma $ is topologically $\mathbb S ^{n-1} \times  (0,1)$ and $\partial \Sigma$ has two connected components $\partial \Sigma _1$ and $\partial \Sigma _2$ homeomorphic to $\mathbb S ^{n-1}$.}
\end{quote}

The Claim A follows as above. The convexity follows from condition \eqref{Eq:PrincipalAnnulus} and the relationship between the principal curvatures and the eigenvalues of the Schouten tensor given in \eqref{lambdakappa}. The boundary components are given by $\partial \Sigma _1 = \phi (\partial \mathbb B ({\bf n}, r))$ and $\partial \Sigma _2 = \phi (\partial \mathbb S ^n _+)$.

\vspace{.3cm}

\begin{quote}
{\bf Claim B: }{\it $\partial \Sigma _1 \subset E(r)$ and $\partial \Sigma _2 \subset E$.}
\end{quote}
\begin{proof}[Proof of Claim B]
That $\partial \Sigma _2 \subset E$ is clear as we did above. It might be that $\partial \Sigma _1 \subset E(r)$ is not that clear. Instead to compute this, we will use a geometric argument. 

Let $i : \mathbb H ^{n+1} \to \mathbb H ^{n+1}$ be the isometry such that $i (E(r)) = E$. Set $\tilde \Sigma = i (\Sigma)$ and $\partial \tilde \Sigma _i = i(\partial \Sigma _i)$, $i=1,2$.

The isometry $i$ induces a conformal diffeomorphism $\Phi :\mathbb S ^n \to \mathbb S ^n $ such that $\Phi (\mathbb S (r)) = \partial \mathbb S ^n _+$. Consider the new metric $\tilde g = \Phi ^* g$, which is defined on 
$$\widetilde{\mathbb{A}(r)} = \Phi (\mathbb A (r))= \mathbb B (\tilde r) \setminus \overline{\mathbb S ^n _+} , \text{ where } \tilde r = r+\pi/2.$$ 

Since \eqref{p2} is conformally invariant, the metric $\tilde g = e^{2 \tilde \rho} g_0$ satisfies 
\begin{equation}\label{p2t}
\left\{\begin{array}{rl}
f(\tilde \lambda(p))= 1,& \tilde \lambda(p)\in\Gamma, \, \, p\in \widetilde{\mathbb{A}(r)},\\ \\

e^{-\tilde \rho}\dfrac{\partial \tilde \rho}{\partial \nu} = 0, & \textrm{ on } \partial \mathbb{S}^n_+ .\\ \\

e^{-\tilde \rho}\dfrac{\partial \tilde \rho}{\partial \nu} + e^{-\tilde \rho}h(\tilde r)= 0, & \textrm{ on } \partial \mathbb{B}(\n ,\tilde r).
\end{array} \right.
\end{equation}

It is not hard to realize that the horospherically concave hypersurface associated to $\tilde g$ is $\tilde \Sigma$ (see \cite{E} for details) with boundary components $\partial \tilde \Sigma _1$ and $\partial \tilde \Sigma _2$.

Thus, as we did above, using \eqref{p2t}, we can check that $\partial \tilde \Sigma _1 \subset E$, which implies that $ \partial \Sigma _1 \subset E(r)$ as claimed. 
\end{proof}

Let us denote by $\mathcal S (r) $ be the connected component of $\mathbb H ^{n+1} \setminus ( E(r) \cup E )$ whose boundary is 
$$ \partial  \mathcal S (r) = E(r) \cup E .$$

\vspace{.3cm}

\begin{quote} 
{\bf Claim C: }{\it $\Sigma \subset \mathcal S (r)$.} 
\end{quote}

The proof of Claim C is as we did in Theorem \ref{thp1}.

\vspace{.3cm}

\begin{quote} 
{\bf Claim D: }{\it $\Sigma$ orthogonal to $E(r)$ (resp. $E$) along $\partial \Sigma_1$ (resp. along $\partial \Sigma _2$).}
\end{quote}

The proof of Claim D is as we did in Theorem \ref{thp1} for each connected component. We also use the isometry $i: \mathbb H ^{n+1} \to \mathbb H ^{n+1}$ of Claim B for $\partial \Sigma _1$.

\vspace{.3cm}

\begin{quote} 
{\bf Claim E: }{\it $\Sigma $ is embedded.}
\end{quote}
\begin{proof}
The proof is as in Claim E in Theorem \ref{thp1} using the foliation $E(s)$ for $s\in [r, \pi/2]$. First, we prove that each boundary component is convex and embedded on the totally geodesic plane that it is contained by \cite{CW}. Then, since they are transversal, there exists $\epsilon >0$ small enough such that $\Sigma _s :=\Sigma \cap E(s)$, $s \in (r, r+2\epsilon) \cup (\pi/2 -2\epsilon ,\pi/2)$ has only one embedded convex component. Therefore, we can close up $\Sigma \setminus \bigcup _{s \in (r, r+\epsilon) \cup (\pi/2 -\epsilon ,\pi/2)} \Sigma _s$ as an embedded convex hypersurface $\tilde \Sigma$. So, $\Sigma$ must be embedded by \cite{CW}.



\end{proof}

\vspace{.3cm}

\begin{quote}
{\bf Claim F: }{\it $\Sigma$ satisfies an elliptic equation $(\mathcal W , \Gamma ^*)$ for a given curvature function depending on its principal curvatures $\kappa (p)= (\kappa _1 (p), \ldots , \kappa _n (p))$; that is, 

\begin{equation}\label{p1W}
\left\{\begin{array}{rl}
\mathcal W (\kappa(p))= 1,& \kappa (p)\in\Gamma ^*, \, \, p\in \Sigma,\\ \\

 \partial \Sigma _1 \perp E(r)  ,\\ \\

\partial \Sigma _2 \perp E ,

\end{array} \right. 
\end{equation}here $\perp$ means that $\Sigma $ is perpendicular to the totally geodesic hyperplanes along its boundary components.}
\end{quote}

The proof of Claim F is exactly the same as in Theorem \ref{thp1}.

\vspace{.3cm}

\begin{quote}
{\bf Claim G: }{\it $\Sigma$ is rotationally symmetric.}
\end{quote}
\begin{proof}
Since $\Sigma $ is orthogonal to $E$ along $\partial \Sigma _2$, one can consider the hyperbolic reflection $R \in {\rm Iso}(\mathbb H ^{n+1})$ across $E$ and double $\Sigma$, that is, $\tilde \Sigma = \Sigma \cup R(\Sigma)$. Now, consider the hyperbolic translation $T_{2r} \in {\rm Iso}(\mathbb H ^{n+1})$ that takes $E(-r)$ into $E(r)$, i.e., $T(E(-r)) = E(r)$. Note that the fixed points at infinity of $T_{2r}$ are the north pole ${\bf n}$ and the south pole ${\bf s}$. 

So, translating $\tilde \Sigma $ using $T_{2r}$, we create a properly embedded hypersurface $\bar \Sigma$ invariant by $T_{2r}$ whose boundary at infinity is $ \partial_{\infty} \bar \Sigma = \{ {\bf s} ,  {\bf n} \} $. Clearly, $\bar \Sigma$ is horospherically concave and satisfies an elliptic equation. 

Take any geodesic $\beta \subset E$ and consider the reflections $R_s : \mathbb H ^{n+1} \to \mathbb H ^{n+1}$ across the totally geodesic hyperplanes orthogonal to $\beta$ at $\beta (s)$, $s \in \mathbb R$. Thus, the Alexandrov Reflection Method applies and we find, for each geodesic $\beta$, a hyperplane of symmetry. Therefore $\bar \Sigma$ is rotationally symmetric. For details see  \cite[Corollary 4.2]{BEQ}. Thus, since $\bar \Sigma$ is rotational symmetric so is $\Sigma$.
\end{proof}

Since $\Sigma$ is rotationally symmetric its horospherical metric $g$ is rotationally symmetric (see \cite{E}). This finishes Theorem \ref{thp2}.

\section{The two dimensional case}

Given a bounded domain $\Omega\subset \mathbb S^2$  the Riemann mapping theorem asserts that  any Riemannian metric $g$ on $\Omega$  is conformal to the round metric  $g_0$ in $\Omega$, say $g=e^{2\rho} g_0$. The problem to find $g$ with constant Gaussian curvature $K$ in $\Omega$ and constant mean curvature $h$ on $\partial\Omega$ is classically referred as the {\it Liouville problem}. In analytical terms, it is equivalent to find a smooth solution to the following problem:

\begin{equation}\label{pde2}
   \left\{ \begin{array}{lrc}
  \De_{g_0} \rho +K=e^{-2\rho} &  \textrm{ in } \Omega, \\ \\
  \frac{\partial \rho}{\partial \nu} = h e^\rho & \textrm{ on } \partial \Omega,
   \end{array}  \right.
  \end{equation}
where $\nu$ is the inward unit vector field to boundary $\partial \Omega.$ 
This problem were solved when $\Omega = \mathbb S^2_+$ by Hang and Wang in \cite{HaWa} and for annuli by Jimenez in \cite{J}.

\smallskip

The Schouten tensor given in equation (\ref{schouten}) is not defined for 2-dimensional metrics. 
However, if we consider a conformal metric $g = e^{2\rho}g_0$ on $\mathbb S^n$, $n>2$, 
we can see that ${\rm Sch}_{g_0} = (1/2)g_0$ and thus
\begin{equation*}\label{schouten2}
{\rm Sch}_g = -\nabla^{2,g_0}\rho + d\rho \otimes d\rho -\frac 1 2 (-1+\| \nabla^{g_0}\rho\|^2_{g_0})g_0.
\end{equation*}

Since this equation makes sense for $n=2$ it can be used as the definition of the Schouten tensor for conformal
metrics on subdomains of $\mathbb S^2$. In fact, it is easy to check that 
$$
{\rm Trace}(g^{-1}{\rm Sch}_g) = -\Delta_g \rho + e^{-2\rho} = 2R(g),
$$ 
where $2R(g)=K(g)$ is the Gaussian curvature function of $g$. 
Thus, if $\lambda_1(p)$ and  $\lambda_2(p)$ denote the eigenvalues of the 
tensor ${\rm Sch}_g $ at $p\in \Omega$, the fully nonlinear form of the  problem (\ref{pde2}) can be written as

\begin{equation}\label{pde3}
   \left\{ \begin{array}{lrc}
  f(\lambda_1(p),\lambda_2(p))=K, & (\lambda_1, \lambda_2)\in\Gamma, \, \, p\in \Omega,\\ \\
  \frac{\partial u}{\partial \nu} = h e^\rho, & \textrm{ on } \partial \Omega,
   \end{array}  \right.
  \end{equation}
where $(f,\Gamma)$ is an elliptic data as defined in Section \ref{intro}.

Applying the techniques we develop in Sections \ref{s3} and \ref{s4} we obtain the following results:

\begin{theorem}\label{thp3} 
Let $\rho \in C^\infty (\overline{\mathbb S ^2 _+})$ be a solution to \eqref{pde3} with
$\Omega = \mathbb {S} ^2 _+$, $K=1$ and $h=c$. 
Then,  $g=e^{2\rho}g_0$ is isometric to the standard round metric $g_0$ on 
$B({\bf n}, r)\subset \mathbb{S}^2$, with $r= {\rm arccot} (c)$.
\end{theorem}

Note that Theorem \ref{thp3} amply generalizes Hang-Wang Theorem \cite{HaWa} to fully nonlinear equations. Moreover, it says that the solution to \eqref{pde3} must be a solution of \eqref{pde2} and hence, in particular, $g$ must be isometric to the standard round metric $g_0$ on  $B({\bf n}, r)\subset\overline{\mathbb{S}^2_+}$.

Now, let us denote by $\overline{\mathbb{A}(r)}$ the annulus $\mathbb {S} ^2 _+\setminus \overline{\mathbb B ({\bf n},r)}$. Then, we can prove:

\begin{theorem}\label{thp4}
Let $\rho \in C^\infty (\overline{\mathbb{A}(r)})$ be a solution to \eqref{pde3} with
$\Omega = \overline{\mathbb{A}(r)}$, 
$K=1$ and $h=0$. 
Then, $g= e^{2\rho}g_0$ is rotationally symmetric metric on $\overline{\mathbb{A}(r)}$. 
\end{theorem}

When $f(\lambda _1 , \lambda _2) = \lambda _1 + \lambda $ we fall into the Liouville Problem \eqref{pde2} studied in \cite{J} in all its generality. In our case, we must focus on the minimality of the boundary, i.e., $h=0$, as in Section \ref{s4}, in contrasts, we study a bigger class of fully nonlinear equations.

\bibliographystyle{amsplain}

\begin{thebibliography}{10}

\bibitem{AE} D.P. Abantos, J.M. Espinar, 
 \newblock{\it In preparation}.



\bibitem{A}
  T. Aubin,
 \newblock {\it \'Equations diff\'erentielles non lin\'eaires et probl\`eme de Yamabe  concernant la courbure scalaire.}
 \newblock  J. Math. Pures Appl. {\bf 55} (1976), 269--296.


\bibitem{BEQ} V. Bonini, J. M. Espinar, J. Qing, 
\newblock{\it Hypersurfaces in in the Hyperbolic Space with support function.} 
\newblock To appear in Adv. in Math. 

\bibitem{CW} M. P. do Carmo, F. W. Warner, \newblock{\it Rigidity and convexity of hypersurfaces in spheres.}  
\newblock J. Differential Geom. {\bf 4}, (1970), 133--144.

\bibitem{CHY} S.Y.A. Chang, Z-C. Han, P. Yang, 
\newblock{\it Classification of singular radial solution to the $\sigma _k$ Yamabe equation on annular domains.} 
J. Differential Equations {\bf 216} (2005), 482--501.

\bibitem{E0} J. F. Escobar,
\newblock{\it Uniqueness theorems on conformal deformation of metrics, Sobolev 
inequalities, and an eigenvalue estimate. }
Comm. Pure Appl. Math. {\bf 43} (1990), no. 7, 857--883.  

\bibitem{E1}
 J. F. Escobar,
 \newblock \emph{Conformal deformation of a Riemannian metric to a 
scalar flat metric with constant mean curvature on the boundary.}
\newblock Ann. of Math. (2) {\bf 136} (1992), 1--50


\bibitem{E2}
 J. F. Escobar,
 \newblock \emph{The Yamabe problem on manifolds with boundary. }
\newblock J. Differential Geom.  {\bf 35} (1992) 21--84.


\bibitem{E} J. M. Espinar, \newblock{\it Invariant conformal metrics on $\mathbb{S}^n$.} Transactions of the A.M.S., {\bf 363} (2011) no. 11, 5649--5662.

\bibitem{EGM} J. M. Espinar, J. A. G\'alvez, P. Mira, \newblock{\it Hypersurfaces in $\mathbb H^{n+1}$ and conformally invariant equations: the generalized Christoffel and Nirenberg problems.} J. Eur. Math. Soc.  {\bf 11} (2009), no. 4, 903--939.

\bibitem{GaMi} J.A. G\'{a}lvez, P. Mira, {\it The Liouville equation in a half-plane.} 
\newblock J. Differential Equations, {\bf 246} (2009), 4173--4187.

\bibitem{HL}
 Z. C. Han and Y.Y. Li,
\newblock \emph{The Yamabe problem on manifolds with boundary: existence and compactness resultus.}
\newblock Duke Math. J. {\bf 99} (1999), 489--542.

\bibitem{HaWa} F. Hang, X. Wang, {\it A new approach to some nonlinear geometric
equations in dimension two.} 
\newblock Calc. Var. Partial Diff. Equations, {\bf 26} (2006), 119--135.

\bibitem{J} A. Jim\'{e}nez, {\it The Liouville equation in an annulus.} 
\newblock J. Nonlinear Anal., {\bf 75} (2012), 2090--2097.

\bibitem{K} N. Korevaar,  
{\it Sphere theorems via Alexandrov for constant Weingarten curvature hypersurfaces. 
Appendix to a note of A. Ros.}
\newblock J. Differential Geom.  {\bf 27} (1988),  221--223.


\bibitem{LL} A. Li, Y.Y. Li,  
\newblock \emph{A fully nonlinear version of the Yamabe problem on manifolds with boundary.} 
\newblock  J. Eur. Math. Soc. (JEMS) {\bf 8} (2006), {\bf no. 2}, 295--316.

\bibitem{M}
 F. C. Marques,
 \newblock \emph{Existence results for the Yamabe problem on manifold with boundary.}
\newblock Indiana Univ. Math. J. {\bf 54} (2005), 1599--1620.

\bibitem{S2}
 R. Schoen,
\newblock \emph{Conformal deformation of a Riemannian metric to constant scalar curvature.}
\newblock J. Diff. Geom. {\bf 20} (1984), 479--495.


\bibitem{T}
N. Trudinger, 
 \newblock \emph{Remarks concerning the conformal deformation of Riemannian structures on compact manifolds.}
\newblock Ann. Scuola Norm. Sup. Pisa (3) {\bf 22} (1968), 265--274.


\bibitem{Y}
H. Yamabe,
\newblock \emph{On a deformation of Riemannian structures on compact manifolds.}
\newblock Osaka Math. J. {\bf 12} (1960), 21--37.

\bibitem{Zha}
L. Zhang,
\newblock \emph{Classification of conformal metrics on $\mathbb R ^2 _+$ with constant Gauss curvature and geodesic curvature on the boundary under various integral finiteness
assumptions.}
\newblock Calc. Var. Partial Diff. Equations 16 (2003), 405–430.


\end{thebibliography}

\end{document}